\documentclass{amsart}

\newtheorem{theorem}{Theorem}[section]
\newtheorem{prop}[theorem]{Proposition}
\newtheorem{cor}[theorem]{Corollary}
\newtheorem{lemma}[theorem]{Lemma}

\theoremstyle{definition}
\newtheorem{definition}[theorem]{Definition}
\newtheorem{example}[theorem]{Example}

\theoremstyle{remark}

\numberwithin{equation}{section}

\DeclareMathOperator{\gl}{GL}
\DeclareMathOperator{\splin}{SL}
\DeclareMathOperator{\psl}{PSL}

\DeclareMathOperator{\id}{Id}

\DeclareMathOperator{\End}{End}

\newcommand{\bmx}{\left(\begin{array}}
\newcommand{\emx}{\end{array}\right)}

\newcommand{\Z}{\mathbb{Z}}
\newcommand{\R}{\mathbb{R}}
\newcommand{\C}{\mathbb{C}}

\input epsf

\begin{document}

\title{Low-Dimensional Unitary Representations of $B_3$}

\author{Imre Tuba}

\address{Department of Mathematics, Mail Code 0112, University of California,
San Diego, 9500 Gilman Dr, La Jolla, CA 92093-0112}

\email{ituba@math.ucsd.edu}

\date{September 5, 1999}


\begin{abstract}
We characterize all simple unitarizable representations of the braid group $B_3$
on complex vector spaces of dimension $d \leq 5$. In particular, we prove that
if $\sigma_1$ and $\sigma_2$ denote the two generating twists of $B_3$, then a
simple representation $\rho:B_3 \to \gl(V)$ (for $\dim V \leq 5$) is
unitarizable if and only if the eigenvalues $\lambda_1, \lambda_2, \ldots,
\lambda_d$ of $\rho(\sigma_1)$ are distinct, satisfy $|\lambda_i|=1$ and
$\mu^{(d)}_{1i} > 0$ for $2 \leq i \leq d$, where the $\mu^{(d)}_{1i}$ are
functions of the eigenvalues, explicitly described in this paper. 
\end{abstract}

\maketitle

\section{Introduction}

Unitary braid representations have been constructed in several ways using the
representation theory of Kac-Moody algebras and quantum groups, see
e.g. \cite{kac}, \cite{squier},
and \cite{wenzl}. Such representations easily lead to representations of
$\psl(2,\Z) = B_3/Z$, where $Z$ is the center of $B_3$, and $\psl(2,\Z) =
\splin(2,\Z) / \{\pm 1\}$, where $\{ \pm 1\}$ is the center
of $\splin(2,\Z)$. We give a complete classification of simple unitary
representations of $B_3$ of dimension $d \leq 5$ in this paper. In particular,
the unitarizability of a braid representation depends only on the the
eigenvalues $\lambda_1, \lambda_2, \ldots, \lambda_d$ of the images the two generating twists of $B_3$. The condition
for unitarizability is a set of linear inequalities in the logarithms of these
eigenvalues. In other words, the representation is unitarizable if and only if
the $(\arg \lambda_1, \arg \lambda_2, \ldots, \arg \lambda_d)$ is a point inside
a polyhedron in $(\R / 2\pi)^d$, where we give the equations of the hyperplanes
that bound this polyhedron. This
classification shows that the approaches mentioned previously do not
produce all possible unitary braid representations. We obtain representations
that seem to be new for $d \geq 3$. As any
unitary representation of $B_n$
restricts to a unitary representation of $B_3$ in an obvious way, these results
may also be useful in classifying such representation of $B_n$.

I would like to thank my advisor, Hans Wenzl for the ideas he contributed to
this paper and Nolan Wallach for his suggestions.

Let $B_3$ be
generated by $\sigma_1$ and $\sigma_2$ with the relation $\sigma_1 \sigma_2
\sigma_1 = \sigma_2 \sigma_1 \sigma_2$. It is well-known that the center of
$B_3$ is generated by $(\sigma_1 \sigma_2)^3$. Let $K$ be any field. If $\rho$
is a simple representation of $B_3$ on a $K$-vector space $V$, then
$\rho(\sigma_1 \sigma_2)^3$ must act on $V$ as a scalar $\delta \in K$. Since
$\sigma_1$ and $\sigma_2$ are conjugates via $\sigma_1 \sigma_2 \sigma_1$,
their images $A=\rho(\sigma_1)$ and $B=\rho(\sigma_2)$ have the same
eigenvalues $\lambda_1, \lambda_2, \ldots, \lambda_d$. We will need the
following two results from \cite{mainbraid}.

\begin{theorem}
\label{mainbraidthm}
\begin{enumerate}
\item
Let $K$ be an algebraically closed field, $V$ a $d$-dimensional $K$-vector
space, and $\lambda_1, \lambda_2, \ldots,\lambda_d \in K-\{0\}$, where $d \leq
5$. There exists a simple representation $\rho:B_3 \to \gl(V)$ such that the
eigenvalues of $A=\rho(\sigma_1)$ satisfy $Q^{(d)}_{rs} \neq 0$ for all $r
\neq s$ where the polynomials $Q^{(d)}_{rs}$ are as follows:
\[ Q^{(2)}_{rs} = -\lambda_r^2 + \lambda_r \lambda_s - \lambda_s^2 \]
\[ Q^{(3)}_{rs} = (\lambda_r^2 + \lambda_s \lambda_k)(\lambda_s^2 +
\lambda_r \lambda_k) \]
with $k \neq r,s$.
\[ Q^{(4)}_{rs} =
-\gamma^{-1}(\lambda_r^2 + \gamma)(\lambda_s^2 + \gamma)(\gamma + \lambda_r
\lambda_k + \lambda_s \lambda_l)(\gamma + \lambda_r \lambda_l + \lambda_s
\lambda_k) \]
with $\gamma = \sqrt{\lambda_1 \cdots \lambda_4}$ and $k,l \neq r,s$.
\[ Q^{(5)}_{rs} = \gamma^{-8}(\gamma^2 + \lambda_r\gamma +
\lambda_r^2)(\gamma^2 + \lambda_s \gamma + \lambda_s^2) \prod_{k \neq r,s}
(\gamma^2+\lambda_r \lambda_k)(\gamma^2+\lambda_s \lambda_k) \]
with $\gamma = \sqrt[5]{\lambda_1 \cdots \lambda_5}$.

\item
A simple representation of $B_3$ of dimension $d \leq 5$ is uniquely
determined up to isomorphism by the eigenvalues of $A=\rho(\sigma_1)$ (for $d
\leq 3$) and $\delta$, where $\rho(\sigma_1 \sigma_2)^3 = \delta\,\id_V$ (for
$d=4,5$).
\end{enumerate}
\end{theorem}

Explicit matrices for $A=\rho(\sigma_1)$ and $B=\rho(\sigma_2)$ are also listed
in \cite{mainbraid}.

The functions $Q^{(d)}_{rs}$ are defined in \cite{mainbraid} by $P_r^{(d)}(B)
P_s^{(d)}(A) P_r^{(d)}(B) = Q^{(d)}_{rs} P_r^{(d)}(B)$, where $P_r^{(d)}(x) =
\prod_{i \neq r} (x-\lambda_i)$. Note that substituting $\lambda_i = e^{2\pi i
t_i}$ and taking logarithms reduces the problem of finding the zeroes of
$Q^{(d)}_{rs}$ to solving a system of linear equations in the $t_i$. (See
Example \ref{exampled=3}.)

\begin{prop}
Let $\rho:B_3 \to \gl(V)$ be a simple representation of dimension $d \leq 5$.
Then the minimal polynomials of $A=\rho(\sigma_1)$ and $B=\rho(\sigma_2)$ are
the same as their characteristic polynomials.
\end{prop}

An immediate consequence of this is

\begin{cor}
\label{diagonalizable}
If $A$ (or $B$) is a diagonalizable matrix, then it has distinct eigenvalues
$\lambda_1, \lambda_2, \ldots, \lambda_d$.
\end{cor}

\begin{proof}
Since $A$ is conjugate to some diagonal matrix $D$, its minimal polynomial is
just $p(x) = \prod (x-d_j)$ where the $d_j$ are the distinct diagonal entries
of $D$. By the previous proposition, $\deg p=d$, hence $D$ must have $d$
distinct diagonal entries. Thus all of the diagonal entries of $D$ are
distinct.
\end{proof}

Since we are interested in unitarizable representations, we will let $K=\C$
and we will require that $|\lambda_i| = 1$. Let $\rho: B_3 \to V$ be a simple
$d$-dimensional representation ($d \leq 5$), and $A=\rho(\sigma_1)$,
$B=\rho(\sigma_2)$. Any unitarizable complex matrix is diagonalizable, so we
can assume that $A$ and $B$ are diagonalizable. So the eigenvalues $\lambda_1,
\lambda_2, \dots, \lambda_d$ are distinct by the last corollary.
Let $\delta$ be the scalar via which $\rho(\sigma_1
\sigma_2)^3$ acts on $V$, that is $(AB)^3 = \delta I$.  Denote the $\C$-algebra
generated by $A$ and $B$ by $\mathcal{B}$. In other words, $\mathcal{B} =
\rho(\C\,B_3$, where $\C\,B_3$ is the group algebra. Note that $\mathcal{B} =
\End(V)$ by simplicity.

The proof proceeds by defining a vector space antihomomorphism
$\imath:\mathcal{B} \to \mathcal{B}$ and proving that it is an algebra
antihomorphism and an involution of $\mathcal{B}$ in section
\ref{involution}. In section \ref{invariance}, we define a sesquilinear
form $\left<.\, ,.\right>$ on the ideal $I=\mathcal{B} e_{B,1}$ that is
invariant under multiplication by $A$ and $B$. We prove that $\left<.\,
,.\right>$ is positive definite if $\mu_{1i}^{(d)} > 0$ for $2 \leq
i \leq d$. In this case, $\rho$ is a unitary representation of $B_3$ on the
$d$-dimensional vector space $I$. We also prove that if $\rho$ is a
unitarizable representation $\mu_{1i}^{(d)} > 0$ for $2 \leq i \leq d$.
In section \ref{examples}, we give some examples of using the positivity of
$\mu_{1i}^{(d)}$.

\section{An involution of the image of $B_3$}
\label{involution}

Let $e_{M,i}$ be the eigenprojection of $M$ to the eigenspace of $\lambda_i$,
where $M\in \{A,B\}$. That is
\[ e_{M,i} = \prod_{j \neq i} \frac{M-\lambda_j}{\lambda_i-\lambda_j} =
\frac{P^{(d)}_i(M)}{\prod_{j \neq i} (\lambda_i-\lambda_j)}. \] 
Note that $e_{A,i}$ and $e_{B,i}$ always exist because the eigenvalues are
distinct. Also $e_{M,i} e_{M,j} = \delta_{ij} e_{M,i}$. Define $\mu^{(d)}_{ij}$
by $e_{B,i} e_{A,j} e_{B,i} = \mu^{(d)}_{ij} e_{B,i}$.  Note that
\[ \mu^{(d)}_{ij} = \frac{Q^{(d)}_{ij}}{\prod_{k \neq i}
(\lambda_i-\lambda_k) \prod_{k \neq j} (\lambda_j-\lambda_k)}.\]

\begin{lemma}
The $\mu^{(d)}_{ij}$ are real numbers. 
\end{lemma}

\begin{proof}
For $i \neq j$, the proof is by direct computation using
$\overline{\lambda_i}=\lambda_i^{-1}$ and $\overline{\gamma}=\gamma^{-1}$. For
example, for $d=5$:

\begin{eqnarray*}
\mu^{(d)}_{ij} & = & \frac{(\gamma^2+\lambda_i\gamma+\lambda_i^2)
(\gamma^2+\lambda_j\gamma+\lambda_j^2) \prod_{k \neq i,j}
(\gamma^2+\lambda_i\lambda_k)(\gamma^2+\lambda_j\lambda_k)}{\gamma^8 \prod_{k
\neq i} (\lambda_i-\lambda_k) \prod_{k \neq j} (\lambda_j-\lambda_k)} \\
& = & \frac{(\gamma \lambda_i^{-1} + 1 + \gamma^{-1} \lambda_i)(\gamma
\lambda_j^{-1} + 1 + \gamma^{-1} \lambda_j)}{(1-\lambda_j
\lambda_i^{-1})(1-\lambda_i \lambda_j^{-1})} \ \frac{\prod_{k \neq i,j}
(\gamma^2+\lambda_i\lambda_k)(\gamma^2+\lambda_j\lambda_k)}{\gamma^6 \prod_{k
\neq i,j} (\lambda_i-\lambda_k)(\lambda_j-\lambda_k)}
\end{eqnarray*}
The first of the two quotients is easily seen to be real. For the second
quotient,
\begin{eqnarray*}
\overline{\left( \frac{\prod_{k \neq i,j}
(\gamma^2+\lambda_i\lambda_k)(\gamma^2+\lambda_j\lambda_k)}{\gamma^6 \prod_{k
\neq i,j} (\lambda_i-\lambda_k)(\lambda_j-\lambda_k)} \right)} & = &
\frac{\prod_{k \neq i,j} (\gamma^{-2}+\lambda_i^{-1}\lambda_k^{-1})
(\gamma^{-2}+\lambda_j^{-1}\lambda_k^{-1})} {\gamma^{-6}\prod_{k
\neq i,j} (\lambda_i^{-1}-\lambda_k^{-1})(\lambda_j^{-1}-\lambda_k^{-1})}
\end{eqnarray*}
Multiply the numerator and the denominator by $\gamma^{12} \lambda_i^3
\lambda_j^3 \prod_{k \neq i,j} \lambda_k^2$ to see that this is still
\[ \frac{\prod_{k \neq i,j}
(\gamma^2+\lambda_i\lambda_k)(\gamma^2+\lambda_j\lambda_k)}{\gamma^6 \prod_{k
\neq i,j} (\lambda_i-\lambda_k)(\lambda_j-\lambda_k)} \]

For the case $i=j$, note that $\sum_{k=1}^d e_{A,k} = I$, so
\begin{eqnarray*}
e_{B,i} & = & e_{B,i} I e_{B,i} \\
& = & e_{B,i} \sum_{k=1}^d e_{A,k} e_{B,i} \\
& = & \sum_{k=1}^d e_{B,i} e_{A,k} e_{B,i} \\
& = & \sum_{k=1}^d \mu^{(d)}_{ik} e_{B,i}
\end{eqnarray*}
Hence $\sum_{k=1}^d \mu^{(d)}_{ik} = 1$, and $\mu^{(d)}_{ii} = 1 - \sum_{k
\neq i} \mu^{(d)}_{ik}$ is real.
\end{proof}

\begin{prop}
$S=\{e_{A,i} e_{B,1} e_{A,j} \mid 1\leq i,j\leq d, \ i \neq j \} \cup \{
e_{A,i} \mid 1 \leq i \leq d \}$ is a basis for the $\C$-vector space $\mathcal{B}$.
\end{prop}

\begin{proof}
Suppose
\[ \sum_{i=1}^d \underset{j \neq i}{\sum_{j=1}^d}
\alpha_{ij} e_{A,i} e_{B,1} e_{A,j} +
\sum_{i=1}^d \alpha_{ii} e_{A,i} = 0 \]
Multiply by $e_{A,i}$ both on the left and on the right. The only term of the
sum that survives is
\[ \alpha_{ii} e_{A,i} = 0 \]
Let $v_i$ be an eigenvector of $A$ corresponding to $\lambda_i$. Then $e_{A,i}
v_i = v_i \neq 0$, so $e_{A,i} \neq 0$. Hence $\alpha_{ii} = 0$.

For $i \neq j$, multiplying by $e_{A,i}$ on the left and by $e_{A,j}$ on the
right shows
\[ \alpha_{ij} e_{A,i} e_{B,1} e_{A,j} = 0 \]
But
\[ e_{B,1} e_{A,i} e_{B,1} e_{A,j} e_{B,1} = (e_{B,1} e_{A,i} e_{B,1})
(e_{B,1} e_{A,j} e_{B,1}) = \mu^{(d)}_{1j} \mu^{(d)}_{1i} e_{B,1} \neq 0 \]
so
$e_{A,i} e_{B,1} e_{A,j} \neq 0$. Hence $\alpha_{ij} = 0$. So $S$ is linearly
independent. It has $d^2$ elements, hence it is a basis of the
$d^2$-dimensional space $\mathcal{B}$.
\end{proof}

Note: if we know $\mu^{(d)}_{ii} \neq 0$ for all $i$, we can use the basis $S'=
\{e_{A,i} e_{B,1} e_{A,j} \mid 1\leq i,j\leq d\}$ instead of $S$. As
$e_{A,i} e_{B,1} e_{A,i} = \mu^{(d)}_{ii} e_{A,i}$, $S'$ is almost the same as
$S$. Since $S'$ is more symmetric than $S$, its use makes the following
computations simpler and the arguments more transparent. In the most
general case however, $\mu^{(d)}_{ii}$ could be $0$.

Define $\imath : \C \to \C$ as the usual complex
conjugation. Extend $\imath$ to $\mathcal{B} \to \mathcal{B}$ by requiring
$\imath$ to be an antilinear map with $\imath(e_{A,i}) = e_{A,i}$ and
$\imath(e_{A,i} e_{B,1} e_{A,j})=e_{A,j} e_{B,1} e_{A,i}$ for $i \neq j$.
Note that $\imath(\mu^{(d)}_{ij})=\mu^{(d)}_{ij}$.

\begin{lemma}
$\imath$ as defined above is an antihomomorphism on the algebra $\mathcal{B}$
and $\imath^2 = \id_{\mathcal{B}}$.
\end{lemma}

\begin{proof}
It is sufficient to prove that $\imath$ acts as an antihomomorphism on the
elements of the basis $S$. $S$ has two different types of elements, therefore
we will have four different cases. Since each can verified directly by a
simple computation, we will show the details for only one:
\begin{enumerate}
\item
\[ \imath(e_{A,i}e_{A,j}) = \imath(e_{A,j})\imath(e_{A,i}) \]

\item
\begin{eqnarray*}
\imath(e_{A,i} (e_{A,j}e_{B,1}e_{A,k})) & = &
\imath(e_{A,j}e_{B,1}e_{A,k})\imath(e_{A,i}) \\
\imath((e_{A,i}e_{B,1}e_{A,j})e_{A,k}) & = &
\imath(e_{A,k})\imath(e_{A,j}e_{B,1}e_{A,k})
\end{eqnarray*}

\item
For $i \neq k$,
\[ \imath((e_{A,i} e_{B,1} e_{A,j})(e_{A,k} e_{B,1} e_{A,l})) = (e_{A,l}
e_{B,1} e_{A,k}) (e_{A,j} e_{B,1} e_{A,i}) \]

\item
\begin{eqnarray*}
\imath((e_{A,i} e_{B,1} e_{A,j})(e_{A,j} e_{B,1} e_{A,k})) & = &
\imath(e_{A,i} (e_{B,1} e_{A,j} e_{B,1}) e_{A,k}) \\ & = & \imath(e_{A,i}
(\mu^{(d)}_{1j} e_{B,1}) e_{A,k}) \\ & = & \overline{\mu^{(d)}_{1j}}
\imath(e_{A,i} e_{B,1} e_{A,k}) \\ & = & \mu^{(d)}_{1j} e_{A,k} e_{B,1}
e_{A,i}
\end{eqnarray*}
Also
\begin{eqnarray*}
\imath(e_{A,j} e_{B,1} e_{A,k}) \imath(e_{A,i} e_{B,1} e_{A,j}) & = & (e_{A,k}
e_{B,1} e_{A,j})(e_{A,j} e_{B,1} e_{A,i}) \\ & = & e_{A,k} (e_{B,1} e_{A,j}
e_{B,1}) e_{A,i} \\ & = & \mu^{(d)}_{1j} e_{A,k} e_{B,1} e_{A,i}
\end{eqnarray*}
\end{enumerate}

That $\imath^2 = \id_{\mathcal{B}}$ follows immediately from the definition.
\end{proof}

\begin{lemma}
\label{eb1}
$\imath(e_{B,1}) = e_{B,1}$.
\end{lemma}

\begin{proof}
First note that $\imath(e_{A,i} e_{B,1} e_{A,i}) = \imath(\mu^{(d)}_{ii}
e_{A,i}) = \mu^{(d)}_{ii} e_{A,i} = e_{A,i} e_{B,1} e_{A,i}$. Multiply
$e_{B,1}$ by $1 = \sum_{i=1}^d e_{A,i}$ on both sides:
\[ e_{B,1}=\left( \sum_{i=1}^d e_{A,i} \right) e_{B,1} \left( \sum_{j=1}^d
e_{A,j} \right) = \sum_{i,j} e_{A,i} e_{B,1} e_{A,j}\] into
\begin{eqnarray*}
\imath(e_{B,1}) & = & \imath \left( \sum_{i=1}^d \sum_{j=1}^d e_{A,i} e_{B,1}
e_{A,j} \right) = \sum_{i=1}^d \sum_{j=1}^d \imath(e_{A,i} e_{B,1} e_{A,j}) \\
& = & \sum_{i=1}^d \sum_{j=1}^d (e_{A,j} e_{B,1} e_{A,i}) = e_{B,1}
\end{eqnarray*}
\end{proof}

\begin{cor}
$\imath(A)=A^{-1}$, and $\imath(I) = I$.
\end{cor}

\begin{proof}
\[ \imath(A) = \imath(\sum_{i=1}^d \lambda_i e_{A,i}) = \sum_{i=1}^d
\overline{\lambda_i} \imath(e_{A,i}) = \sum_{i=1}^d \lambda_i^{-1} e_{A,i} =
A^{-1} \] Similarly,
\[ \imath(I) = \imath(\sum_{i=1}^d e_{A,i}) = \sum_{i=1}^d \imath(e_{A,i}) =
\sum_{i=1}^d e_{A,i} = I \]
\end{proof}

\begin{lemma}
$\imath(B)=B^{-1}$.
\end{lemma}

\begin{proof}
Note that $A^{-1} \imath(B) A^{-1} = \imath(A) \imath(B) \imath(A) =
\imath(ABA) = \imath(BAB) = \imath(B) A^{-1} \imath(B)$. That is $A^{-1}$ and
$\imath(B)$ satisfy the braid relation. So the group homomorphism $\rho':B_3
\to \gl(V)$ defined by $\rho'(\sigma_1)=A^{-1}$ and
$\rho'(\sigma_2)=\imath(B)$ is another representation of $B_3$ on $V$. Once
again, the braid relation implies that $A^{-1}$ and $\imath(B)$ are
conjugates. Hence they have the same eigenvalues, namely $\lambda_1^{-1},
\lambda_2^{-1}, \ldots, \lambda_d^{-1}$.

But $\imath: \mathcal{B} \to \mathcal{B}$ only permutes the basis $S$ of
$\mathcal{B}=\End(V)$.
Hence $\imath(\mathcal{B}) = \imath(\End(V)) = \End(V)$ and $A^{-1}$ and
$\imath(B)$ generate the algebra $\End(V)$. That is $\rho'$ is also a simple
representation of $B_3$

Now, $(A^{-1} \imath(B))^3 = \imath(BA)^3 = \imath(AB)^3 = \imath(\delta I) =
\overline{\delta} = \delta^{-1} I$ (recall $|\delta|=1$). By Corollary
\ref{mainbraidthm}, the
eigenvalues $\lambda_1^{-1}, \lambda_2^{-1}, \ldots, \lambda_d^{-1}$ (if
d=$2,3$) or the eigenvalues together with $\delta$ (if $d=4,5$) uniquely
determine a simple representation of $B_3$ on $V$ up to isomorphism.

But we already know such a representation, namely $\sigma_1 \mapsto A^{-1}$
and $\sigma_2 \mapsto B^{-1}$. Hence there exists $M \in \gl(V)$ such that
$A^{-1}=MA^{-1}M^{-1}$ and $\imath(B)=MB^{-1}M^{-1}$. Then $M$ is in the
centralizer of $A$.

\begin{eqnarray*}
Me_{B,1}M^{-1} & = & M \left( \prod_{i=2}^d
\frac{B-\lambda_i}{\lambda_1-\lambda_i} \right) M^{-1} \\ & = & \prod_{i=2}^d
\frac{MBM^{-1}-\lambda_i}{\lambda_1-\lambda_i} \\ & = & \prod_{i=2}^d
\frac{\imath(B^{-1})-\lambda_i}{\lambda_1-\lambda_i} \\ & = & \prod_{i=2}^d
\imath \left(\frac{B^{-1}-\lambda_i^{-1}}{\lambda_1^{-1}-\lambda_i^{-1}}
\right) \\ & = & \imath \left( \prod_{i=2}^d
\frac{B^{-1}-\lambda_i^{-1}}{\lambda_1^{-1}-\lambda_i^{-1}} \right)
\end{eqnarray*}
Call the quantity in parentheses $\phi$. Note that $\phi$ is the
eigenprojection to the subspace spanned by the eigenvector $w_1$ of $B^{-1}$
with eigenvalue $\lambda_1^{-1}$. But the eigenvectors $w_1, w_2, \ldots,
w_d$ of $B^{-1}$ are also eigenvectors of $B$ and span $V$ (the eigenvalues
are distinct). Hence $\phi(w_1) = w_1 = e_{B,1} w_1$ and $\phi(w_i) = 0 =
e_{B,1} w_i$ for $i \geq 2$. That is $\phi=e_{B,1}$ as their action on the
basis $\{ w_1, w_2, \ldots, w_d \}$ is identical. Then
Lemma \ref{eb1} shows $\imath(Me_{B,1}M^{-1}) = \imath(\phi) = \imath(e_{B,1})
= e_{B,1}$.

Hence conjugation by $M$ is a $\mathcal{B}$-algebra isomorphism that fixes $A$
and $e_{B,1}$. But $A$ and $e_{B,1}$ generate the basis $S$ of
$\mathcal{B}$, hence they generate the algebra $\mathcal{B}$. So conjugation
by $M$ must fix every element of $\mathcal{B}$. In particular, $\imath(B) =
MB^{-1}M^{-1} = B^{-1}$.
\end{proof}

\section{An invariant inner-product}
\label{invariance}

Let $\mathcal{B}$ act on the left algebra ideal $\mathcal{B}e_{B,1}$. Note
that $\mathcal{B}e_{B,1}$ is a $d$-dimensional $\C$-vector space, as $e_{B,1}$
is an idempotent of rank $1$.

\begin{definition}
Define the form $\left<.\, ,.\right>$ on $\mathcal{B}e_{B,1}$ by
$\left<ae_{B,1},be_{B,1} \right> e_{B,1} = \imath(be_{B,1}) ae_{B,1} = e_{B,1}
\imath(b) ae_{B,1}$ for $ae_{B,1}, be_{B,1} \in \mathcal{B}e_{B,1}$.
\end{definition}

It is easy to verify that $\left<.\, ,.\right>$ is a sesquilinear form on the
$\C$-vector space $\mathcal{B}e_{B,1}$. Since $\imath(A)=A^{-1}$ and
$\imath(B)=B^{-1}$, this form is clearly invariant under the action by $A$ and
$B$, hence $\rho(B_3)$.

\begin{lemma}
\label{basis-t}
$T=\{e_{A,i}e_{B,1} \mid 2 \leq i \leq d \} \cup \{ABAe_{B,1}\}$ is a basis for
left algebra ideal $\mathcal{B}e_{B,1}$ considered as a $\C$-vector space.
\end{lemma}

\begin{proof}
Suppose
\[ \alpha_1 ABAe_{B,1} + \sum_{i=2}^d \alpha_i e_{A,i}e_{B,1} = 0 \]

Note that $(e_{A,i} ABA e_{B,1}) (ABA)^{-1} = e_{A,i} e_{A,1} =
\delta_{1i}$. Since $(ABA)^{-1}$ is invertible $e_{A,i} ABA e_{B,1} = 0$ if
and only if $i \geq 2$.

Multiply by $e_{A,1}$ on the left. Then $\alpha_1 e_{A,1} ABA e_{B,1} = 0$
But $e_{A,1} ABA e_{B,1} \neq 0$, so $\alpha_1 = 0$.

Now, multiply by $e_{A,i}$ ($i \geq 2$) on the left. Then $\alpha_i e_{A,i}
e_{B,1} = 0$. We know $e_{B,1} e_{A,i} e_{B,1} = \mu^{(d)}_{1i} e_{B,1} \neq
0$ by simplicity, so $e_{A,i} e_{B,1} \neq 0$ and $\alpha_i = 0$.

Hence $T$ is a linearly independent set, and we can conclude that it is a basis
of the $d$-dimensional vector space $\mathcal{B} e_{B,1}$.
\end{proof}

Note: if we know $e_{A,1} e_{B,1} \neq 0$, we can use the more symmetric basis
$T' = \{ e_{A,i} e_{B,1} \mid 1 \leq i \leq d \}$ to simplify this and some of
the following computations. Unfortunately, $e_{A,1} e_{B,1}$ could in general
be $0$. In particular, if $\mu^{(d)}_{11} = 0$, then $e_{A,1} e_{B,1} = 0$ too.

\begin{theorem}
The braid representation $\mathcal{B}$ is unitarizable if and only if
$\mu^{(d)}_{1i}>0$ for all $2 \leq i \leq d$.
\end{theorem}

\begin{proof}
Suppose $\mu^{(d)}_{1i}>0$ for all $2 \leq i \leq d$. Consider the action of
$\mathcal{B}$ on $\mathcal{B}e_{B,1}$. The sesquilinear form defined above is
invariant under the action of $\rho(B_3)$. So it is sufficient to show that it
is an inner product. That is we need to prove is that it is positive
definite. On the basis $T$:
\begin{eqnarray*}
\left< e_{A,i}e_{B,1}, e_{A,i}e_{B,1} \right> e_{B,1} & = & e_{B,1}
\imath(e_{A,i}) e_{A,i} e_{B,1} = e_{B,1} e_{A,i} e_{A,i} e_{B,1} \\
& = & e_{B,1} e_{A,i} e_{B,1} = \mu^{(d)}_{1i} e_{B,1} \\
\left< ABAe_{B,1}, ABAe_{B,1} \right> e_{B,1} & = & \left< e_{B,1}, e_{B,1}
\right> e_{B,1} = e_{B,1}e_{B,1} = e_{B,1}
\end{eqnarray*}
Hence $\left< e_{A,i}e_{B,1}, e_{A,i}e_{B,1} \right> = \mu^{(d)}_{1i} > 0$ for
$i \geq 2$ by assumption, and $\left< ABAe_{B,1},\, ABAe_{B,1} \right>=1$.
We claim that $T$ is orthogonal with respect to $\left<.,\, ,\right>$. Let $i,j
\neq 1$ and $i \neq j$:
\begin{eqnarray*}
\left< e_{A,i}e_{B,1}, e_{A,j}e_{B,1} \right> e_{B,1} & = & e_{B,1}
\imath(e_{A,i}) e_{A,j} e_{B,1} = e_{B,1} e_{A,i} e_{A,j} e_{B,1} = 0 \\
\left< ABAe_{B,1}, e_{A,i}e_{B,1} \right> e_{B,1} & = & e_{B,1}
\imath(e_{A,i}) ABA e_{B,1} = e_{B,1} e_{A,i} ABA e_{B,1} = 0
\end{eqnarray*}
We used $e_{A,i} ABA e_{B,1} = 0$ in the last computation just like in Lemma
\ref{basis-t}.

Hence $\left<.,\, .\right>$ is a positive definite form. Then
$\mathcal{B}e_{B,1}$ is a $\C$-vector space with inner product $\left<.\,
,.\right>$ and the action of $\rho(B_3)$ on this space is unitary.

Conversely, suppose $\mathcal{B}$ is unitarizable. So there exists $V$ a $\C$
vector space with inner product $\left<.\, ,.\right>$ and $\rho: B_3 \to
\gl(V)$ such that $A=\rho(\sigma_1)$ and $B=\rho(\sigma_2)$ act as unitary
operators on $V$. Let $^*$ be the transpose induced by $\left<.\,
,.\right>$. We know $A^*=A^{-1}$ and $B^*=B^{-1}$. Let $v \in V$ be an
eigenvector of $B$ with eigenvalue $\lambda_1$. Then $e_{B,1}v=v$ and
\begin{eqnarray*}
0 \leq \left<e_{A,i}e_{B,1}v,e_{A,i}e_{B,1}v \right> & = & \left< v,
e_{B,1}^*e_{A,i}^*e_{A,i}e_{B,1}v \right> \\ & = & \left<
v,e_{B,1}e_{A,i}e_{B,1}v \right> = \left< v, \mu^{(d)}_{1i} e_{B,1} v \right>
= \mu^{(d)}_{1i} \left< v,v \right>
\end{eqnarray*}
Hence $\mu^{(d)}_{1i} \geq 0$. We know $\mu^{(d)}_{1i} \neq 0$ for $i\geq 2$ by
simplicity, so $\mu^{(d)}_{1i}>0$ in this case.
\end{proof}

\section{Examples}
\label{examples}

\begin{example}
$d=2$
\end{example}
\begin{eqnarray*}
\mu^{(2)}_{12} & = & \frac{-\lambda_1^2 + \lambda_1 \lambda_2 -
\lambda_2^2}{(\lambda_1-\lambda_2)(\lambda_2-\lambda_1)}\\
& = & \frac{\lambda_1^2 - \lambda_1 \lambda_2 +
\lambda_2^2}{(\lambda_1-\lambda_2)^2} \\
& = & 1 + \frac{\lambda_1\lambda_2}{(\lambda_1-\lambda_2)^2}\\
& = & 1 - \frac{1}{(\lambda_1/\lambda_2-1)(\lambda_2/\lambda_1-1)}\\
& = & 1 - \left| \frac{\lambda_1}{\lambda_2} - 1 \right|^{-2} > 0
\end{eqnarray*}
That is
\[ \left| \frac{\lambda_1}{\lambda_2} - 1 \right| > 1 \]
or $\lambda_1/\lambda_2 = e^{it}$ for $\pi/3 < t <5\pi/3$.

\begin{example}
\label{exampled=3}
$d=3$
\end{example}
\begin{eqnarray*}
\mu^{(3)}_{12} & = & \frac{(\lambda_1^2 + \lambda_2 \lambda_3)(\lambda_2^2 +
\lambda_1 \lambda_3)}
{(\lambda_1-\lambda_2)(\lambda_1-\lambda_3)(\lambda_2-\lambda_1)(\lambda_2-\lambda_3)} \\
& = & \frac{\left(1 + \frac{\lambda_3}{\lambda_1}\frac{\lambda_2}{\lambda_1}\right)
\left(\frac{\lambda_2}{\lambda_1}+\frac{\lambda_3}{\lambda_2}\right)}
{\left(1-\frac{\lambda_2}{\lambda_1}\right)
\left(1-\frac{\lambda_1}{\lambda_2}\right)
\left(1-\frac{\lambda_3}{\lambda_1}\right)
\left(\frac{\lambda_2}{\lambda_1}-\frac{\lambda_3}{\lambda_1}\right)} \\
\mu^{(3)}_{13} & = & \frac{(\lambda_1^2 + \lambda_2 \lambda_3)(\lambda_3^2 +
\lambda_1 \lambda_2)}
{(\lambda_1-\lambda_2)(\lambda_1-\lambda_3)(\lambda_3-\lambda_1)(\lambda_3-\lambda_2)} \\
& = &
\frac{\left(1 + \frac{\lambda_2}{\lambda_1}\frac{\lambda_3}{\lambda_1}\right)
\left(\frac{\lambda_3}{\lambda_1}+\frac{\lambda_2}{\lambda_3}\right)}
{\left(1-\frac{\lambda_3}{\lambda_1}\right)
\left(1-\frac{\lambda_1}{\lambda_3}\right)
\left(1-\frac{\lambda_2}{\lambda_1}\right)
\left(\frac{\lambda_3}{\lambda_1}-\frac{\lambda_2}{\lambda_1}\right)}
\end{eqnarray*}
Let $\omega_2 = \lambda_2/\lambda_1$ and $\omega_3 =
\lambda_3/\lambda_1$. Then
\begin{eqnarray*}
\mu^{(3)}_{12} & = & 
\frac{\left(1 + \omega_3
\omega_2\right)\left(\omega_2+\omega_3\omega_2^{-1}
\right)}
{\left|1-\omega_2\right|^2 \left(1-\omega_3\right)
\left(\omega_2-\omega_3\right)} \\
\mu^{(3)}_{13} & = & 
\frac{\left(1 + \omega_2 \omega_3\right)\left(\omega_3+\omega_2\omega_3^{-1}
\right)}
{\left|1-\omega_3\right|^2 \left(1-\omega_2\right)
\left(\omega_3-\omega_2\right)}
\end{eqnarray*}

Let $e^{2\pi t_2} = \omega_2$ and $e^{2\pi t_3} = \omega_3$. So we are looking
for $(t_2,t_3) \in [0,1)^2$ such that both $\mu^{(3)}_{12} > 0$ and
$\mu^{(3)}_{13} > 0$. $\mu^{(3)}_{12}$ and $\mu^{(3)}_{13}$ can change signs
at
\begin{eqnarray*}
\omega_2 \omega_3 & = & -1 \\
\omega_3 \omega_2^{-1} & = & - \omega_2  \\
\omega_2 \omega_3^{-1} & = & - \omega_3  \\
w_2 & = & 1 \\
w_3 & = & 1 \\
w_2 & = & w_3
\end{eqnarray*}
These equations can be transformed into linear equations in $t_2$ and $t_3$ by
taking logs:
\begin{eqnarray*}
t_2 + t_3 & = & \frac{1}{2} \\
t_3 & = & 2t_2 + \frac{1}{2} \\
t_2 & = & 2t_3 + \frac{1}{2} \\
t_2 & = & 0 \\
t_3 & = & 0 \\
t_2 & = & t_3
\end{eqnarray*}
Of course, the above equations are all understood mod $1$.

Computation by Maple shows that $\mu^{(3)}_{12} > 0$ and $\mu^{(3)}_{13} > 0$
in the open set colored black on the plot below. The grey regions are those
where one of $\mu^{(3)}_{12}$ and $\mu^{(3)}_{13}$ is positive and the other is
negative. The line $t_2=t_3$ corresponds to $\lambda_2=\lambda_3$, in which
case the representation cannot be unitarizable.

\[ \epsffile{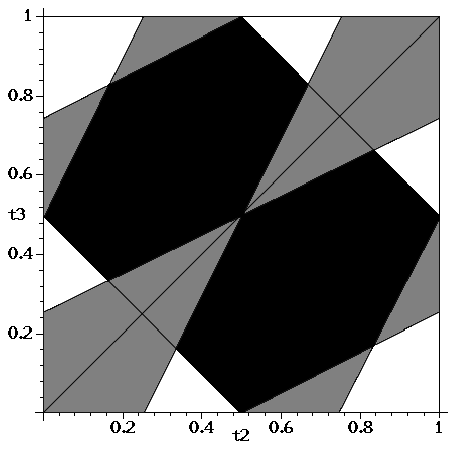} \]

\end{document}